\documentclass{slides}

\usepackage{amsmath,amssymb}
\usepackage{latexsym}

\newcommand{\GL}{GL_n(\mathbb{R})}
\newcommand{\BZ}{\mathbb{Z}}
\newcommand{\BQ}{\mathbb{Q}}

\newcommand\g{\gamma}

\newcommand\G{\Gamma}

\newcommand\la{\lambda}

\newcommand\s{\sigma}

\newcommand\lrt{\longrightarrow}

\newcommand\sd{\,\,{\vartriangleright}\kern -1.0ex{<}\,}

\newcommand\CP{{\mathcal P}_g}

\newcommand\BD{\mathbb D}
\newcommand\BH{\mathbb H}

\begin{document}
\title{\Huge \bf Lecture on Langlands Functoriality Conjecture }
\author{{\Large \bf Jae-Hyun Yang}\\ \vskip 0.5cm {\large \bf Inha University} \\
jhyang@inha.ac.kr}
\vskip 2cm

\
\date{Department of Mathematics\\ Kyoto University\\
 Kyoto, Japan
\\ June 30 (Tue), 2009}


\def\lrt{\longrightarrow}
\def\CP{{\bf{\cal P}}_n}
\def\CPR{{\cal P}_n\times {\Bbb R}^{(m,n)}}
\def\SP{S{\cal P}_n}
\def\GL{GL_n(\Bbb R)}
\def\Rmn{{\Bbb R}^{(m,n)}}
\def\BR{\Bbb R}
\def\BC{\Bbb C}
\def\g{\gamma}
\def\G{\Gamma}
\def\O{\Omega}
\def\o{\omega}
\def\la{\lambda}
\def\BD{\Bbb D}
\def\BH{\Bbb H}
\def\Cmn{{\Bbb C}^{(m,n)}}
\def\PZ{ {{\partial}\over {\partial Z}} }
\def\PW{ {{\partial}\over {\partial W}} }
\def\PZB{ {{\partial}\over {\partial{\overline Z}}} }
\def\PWB{ {{\partial}\over {\partial{\overline W}}} }
\def\PO{ {{\partial}\over {\partial \Omega}} }
\def\PE{ {{\partial}\over {\partial \eta}} }
\def\POB{ {{\partial}\over {\partial{\overline \Omega}}} }
\def\PEB{ {{\partial}\over {\partial{\overline \eta}}} }
\def\PX{ {{\partial}\over{\partial X}} }
\def\PY{ {{\partial}\over {\partial Y}} }
\def\PU{ {{\partial}\over{\partial U}} }
\def\PV{ {{\partial}\over{\partial V}} }
\def\HC{\mathbf{H}_n\times \BC^{(m,n)} }
\def\A{\mathbb A}
\newcommand\LH{{}^LH}
\newcommand\LG{{}^LG}
\newcommand\GA{G({\mathbb A})}

\maketitle

\newpage
\begin{slide}
\begin{center}
{\large \bf 1. Langlands Functoriality Conjecture (briefly LFC)}
\end{center}

\vskip 0.1cm
$F$ = number field\\
 $G$ : connected,
quasi-split reductive /$F$. \par
For each
place $v$ of $F$, \par
 $F_v$ : Completion of $F$ at $v$,\\
${\mathfrak o}_v$ : Ring of integers of $F_v$,\\
${\mathfrak
p}_v$ : Maximal ideal of ${\mathfrak o}_v$,\\
 $q_v$ : Order of residue field $k_v={\mathfrak o}_v /{\mathfrak p}_v$\\
 ${\mathbb A}={\mathbb A}_F$ : Ring of adeles of
$F$.\par For each place $v$ of $F$, we let $G_v=G(F_v).$
\end{slide}

\newpage
\begin{slide}
\vskip 0.1cm Let $$\Omega=\otimes_v \Omega_v\,:\ {\rm Cusp\ Auto\ Repn\ of}\ \GA$$
Let $S$ be a finite set
of places including all archimedean ones such that both $\Omega_v$
and $G_v$ are unramified for any place $v\notin S$. Then for each
$v\notin S,\ \Omega_v$ determines uniquely a semi-simple conjugacy
class
$$c(\Omega_v)\subset {}^LG_v$$ in the $L$-group ${}^LG_v$ of $G_v$ as a group
defined over $F_v$. We note that there exists a
natural homomorphism $$\xi_v:{}^LG_v\lrt {}^LG.$$ For a finite
dimensional repn $r$ of $\LG$, putting $r_v=r\circ
\xi_v$, the {\bf local Langlands $L$-function} $L(s,\Omega_v,r_v)$
associated to $\Omega_v$ and $r_v$ is defined to be
\begin{equation*}
L(s,\Omega_v,r_v)=\det \Big( I-r_v\big(c(\Omega_v)\big)
q_v^{-s}\Big)^{-1}.
\end{equation*}
\end{slide}

\newpage
\begin{slide}

We set
\begin{equation*}
L_S(s,\Omega,r)=\prod_{v\notin S}L(s,\Omega_v,r_v).
\end{equation*}

Langlands [1970] proved that $L_S(s,\Omega,r)$ converges
absolutely for sufficiently large ${\rm Re} (s)> 0$ and defines a
holomorphic function there. Furthermore he proposed the following
question. \vskip 0.1cm \noindent {\bf Conjecture
(Langlands\,[1970]).} $L_S(s,\Omega,r)$ has a meromorphic
continuation to the whole complex plane and satisfies a standard
functional equation. \vskip 0.21cm

F. Shahidi (cf.\,[1988]\,[1990]) gave a partial answer
to the above conjecture using the so-called {\it Langlands-Shahidi
method}.
\end{slide}

\newpage
\begin{slide}
{\bf Special Case $G=GL(n)$:} \par
[1] (Godement-Jacquet[1972], Jacquet[1979])

Let $n\geq 1$ be a positive integer. Assume
$\pi=\otimes_v \pi_v$ be a nontrivial cuspidal representation of $GL(n,{\mathbb A}_F)$.
Then $L(s,\pi):=L(s,\pi, {\mathrm{st}})$ is {\it entire}.
Moreover, for any finite set $S$ of places of $F$, the incomplete $L$-functions
\begin{equation*}
L_S(s,\pi)=\prod_{v\notin S}L(s,\pi_v)
\end{equation*}
is holomorphic and non-zero in ${\mathrm {Re}}(s) >w+1$ if $\pi$ has weight $w$. And we have a
functional equation
\begin{equation*}
L(w+1-s,\pi^{\vee})=\,\varepsilon(s,\pi)\,L(s,\pi)
\end{equation*}
with
$$ \varepsilon(s,\pi)=\,W(\pi)\,N_\pi^{(w+1)/2-s}.$$
Here $\pi^{\vee}$ denotes the contragredient of $\pi$, $N_{\pi}$ denotes the norm of the conductor ${\mathcal N}_\pi$ of $\pi$, and $W(\pi)$ is
the root number of $\pi.$ When $n=1$, this result is due to {\bf Hecke}.

\end{slide}


\newpage
\begin{slide}
[2] (Jacquet-Shalika[1981])

Given any $m$-tuple of cuspidal representations $\pi,\cdots,\pi_m$ of $GL(n_1,{\mathbb A}_F),\cdots,GL(n_m,{\mathbb A}_F)$ respectively,
there exists an irreducible automorphic representation
$$ \pi_1\boxplus\pi_2 \boxplus\cdots \boxplus \pi_m \ \ ({\mathrm {{\bf isobaric}\ auto\ reprn}})$$
of $GL(n,{\mathbb A}_F),\ n=\sum_{j=1}^mn_j$, which is unique up to equivalence such that for any finite set $S$ of places,
\begin{equation*}
L_S(s,\boxplus_{j=1}^m \pi_j)=\,\prod_{j=1}^m L_S(s,\pi_j).
\end{equation*}
Here $\boxplus$ denotes the Langlands sum that comes from his theory of Eisenstein series.

[3] (Jacquet-Shalika[1981])

Let $\pi$ be an automorphic reprn of $GL(n,{\mathbb A}_F)$. Then the order of pole at $s=1$ of
$L_S(s,\pi\times \pi^{\vee})$ is {\bf $1$} if and only if $\pi$ is {\bf cuspidal}.

\end{slide}


\newpage
\begin{slide}
\noindent{\large{\bf Local Langlands Conjecture}} \\
{\bf (briefly,} {\large{\bf LLC})}

Let $G$ be a connected reductive group over a local field $k$.
Let ${\mathcal G}_k(G)$ be the set of all admissible homomorphisms
$\phi:W_k'\lrt \,{}^LG$ modulo inner automorphisms by elements of ${}^LG^0$, and let
$\prod (G(k))$ be the set of all equivalence classes of irreducible admissible representations of $G(k).$
Then there is a surjective map
$$\prod (G(k))\lrt {\mathcal G}_k(G)$$
with finite fibres which partitions $\prod (G(k))$ into disjoint finite sets $\prod_\phi (G(k))$,
called $L$-{\bf packets} satisfying suitable properties. Here
$$W_k'=\,W_k\ltimes {\mathbb G}_a\quad ({\textrm{semidirect product}})$$
denotes the Weil-Deligne group ($W_k$ is the Weil group of $k$).

\end{slide}


\newpage
\begin{slide}

\noindent {\bf Remark:}
(a) If $k$ is
archimedean, i.e., $k=\BR$ or $\BC$, [LLC] was solved by Langlands
(1973).
\par
(b) In case $k$ is non-archimedean, Kazhdan and Lusztig (1987)
had shown how to parametrize those irreducible admissible
representations of $G(k)$ having an Iwahori fixed vector in terms
of admissible homomorphisms of $W_k'$.
\par
(c) For a local field $k$ of positive characteristic $p>0,$ [LLC]
was established by Laumon, Rapoport and Stuhler (1993).
\par
(d) In case $G=GL(n)$ for a non-archimedean local field $k$, [LLC]
was established by Harris and Taylor (2001), and by Henniart (2000).
In both cases, the correspondence was established at
the level of a correspondence between irreducible Galois
representations and supercuspidal representations.

\end{slide}


\newpage
\begin{slide}

(e) Let $k$ be a a non-archimedean local field of characteristic
$0$ and let $G=SO(2n+1)$ the split special orthogonal group over
$k$. In this case, Jiang and Soudry (2003) gave a
parametrization of {\it generic} {\bf supercuspidal} representations of
$SO(2n+1)$ in terms of admissible homomorphisms of $W_k'.$ More
precisely, there is a unique bijection of the set of conjugacy
classes of all admissible, completely reducible,
multiplicity-free, symplectic complex representations
$\phi:W_k'\lrt {}^LSO(2n+1)=Sp(2n,\BC)$ onto the set of all
equivalence classes of irreducible generic {\bf supercuspidal}
representations of $SO(2n+1,k).$ \\
We recall that $(\pi,V)$ is said to be
{\it supercuspidal} if the matrix coefficient
$$f_{v,v^*}(g):=(\pi(g)v,v^*),\quad g\in G$$
has compact support on $G$ modulo the center $Z$ of $G$, where $v\in V$ and
$v^*\in V^*$\,(the dual space of $V$) or equivalently if the Jacquet module $V_N=0$
for any unipotent radical $N$.

\end{slide}


\newpage
\begin{slide}

\vskip 0.3cm For $\pi\in \prod_\phi(G)$ with $\phi\in {\mathcal G}_k(G),$ if $r$
is a finite dimensional complex representation of $\LG$, we define
the $L$- and $\varepsilon$-factors
\begin{equation*}
L(s,\pi,r)=L(s,r\circ \phi)
\end{equation*}
and
$$\varepsilon(s,\pi,r,\psi)=\varepsilon(s,r\circ \phi,\psi),$$
\noindent where $L(s,r\circ \phi)$ is the Artin-Weil $L$-function.

\vskip 0.3cm\noindent {\bf Remark:} The representations in the
$L$-packet $\prod_\phi$ are parametrized by the component group
$$C_\phi:=S_\phi/ Z_L S_\phi^0,$$
where $S_\phi$ is the centralizer of the image of $\phi$ in $\LG,\
S_\phi^0$ is the identity component of $S_\phi$, and $Z_L$ is the
center of $\LG$.

\end{slide}


\newpage
\begin{slide}
\noindent {\huge\bf LFC.} Let $F$ be a number field, and let $H$ and $G$ be connected reductive groups over $F$ such that $G$
is quasi-split over $F$. Suppose $$\sigma: {}^LH \lrt {}^LG$$ is an $L$-homomorphism.
Then for any
automorphic repn $\pi=\otimes_v \pi_v$ of $H(\A)$, there
exists an automorphic repn $\Pi=\otimes_v \Pi_v$ of
$G(\A)$ such that
\begin{equation}
c(\Pi_v)=\s(c(\pi_v)),\quad v\notin S(\pi)\cup S(\Pi),
\end{equation}
where
$$S(\pi)(resp.\ S(\Pi))=\big\{ v\,|\ \pi_v \,(resp.\,\Pi_v)\ \textrm{is\ ramified}\,\big\}.$$

We note that (1) is equivalent
to the condition
\begin{equation*}
L_S(s,\Pi,r)=L_S(s,\pi,r\circ \s),\quad S=S(\pi)\cup S(\Pi)
\end{equation*}
for every finite dim complex repn $r$ of $\LG$.
\end{slide}

\newpage
\begin{slide}
\begin{center}
{\large \bf 2. Examples of Functoriality}
\end{center}

\noindent{\large\bf Example 1.} Let $H=\{1\}$ and $G=GL(n)$. Then $$\LH=\Gamma_F=Gal\left( {\bar F}/F\right),\ \
\LG=GL(n,\BC)\rtimes \Gamma_F.$$ We set
$${\bar \sigma}:Gal(K/F)\lrt GL(n,\BC)\lrt PGL(n,\BC).$$
(1) If $n=1$, $F^{\times}\cong W_F^{\rm ab}$. This is the so-called {\bf Artin Reciprocity Law} . In this case,
 ${\rm Im}\,{\bar\sigma}$ is cyclic or dihedral (E. Artin, E. Hecke)
\par\noindent (2) The case $n=2$. The image ${\rm Im}\,{\bar\sigma}$ is classified as follows:\\
(2-a) ${\rm Im}\,{\bar\sigma}$ is cyclic or dihedral (Artin).\\
(2-b) ${\rm Im}\,{\bar\sigma}=A_4$ is tetrahedral (R. Langlands; {\bf Base Change} for $GL(2)$ [1980]).\\
(2-c) ${\rm Im}\,{\bar\sigma}=S_4$ is octahedral (J. Tunnell [1981]).
\end{slide}

\begin{slide}
\newpage
(2-d) ${\rm Im}\,{\bar\sigma}=A_5$ is icosahedral (R. Taylor et al.[2003] and C. Khare [2005]).\\

(3) The case $Gal(K/F)$ is {\bf nilpotent} and $K/F$ cyclic of prime degree for
$n\geq 3$ (L. Clozel and J. Arthur : {\bf BC} for $GL(n)$ [1989]).

\smallskip

\noindent {\large\bf Example 2.} (Cogdell, Kim, P.-S, Shahidi [2004])\\
(1) $H=Sp(2n),\ G=GL(2n+1)$ and $\sigma$ is an embedding.\\
(2) $H=SO(2n+1),\ G=GL(2n)$ and $\sigma$ is an embedding.\\
(c) $H=SO(2n),\ G=GL(2n)$ and $\sigma$ is an embedding.\\
All the above classical groups are of the {\it split} form.
For an irreducible, {\bf generic} cuspidal repn $\pi$ of $H(\A)$, there is a {\it weak} functorial lift of $\pi$.
The proof is based on {\bf Converse Theorem} of Cogdell and Piateski-Shapiro.
\end{slide}

\begin{slide}
\newpage

Let $\pi=\,\otimes_v\pi_v$ be an irreducible admissible reprn of $GL(n,{\mathbb A})$ and $\tau=\otimes_v\tau_n$
be a cuspidal reprn of $GL(m,{\mathbb A})$ with $m<n.$ Let $\psi$ be a fixed nontrivial continuous additive
character of ${\mathbb A}$ which is trivial on $F$. We define formally
$$L(s,\pi\times\tau)=\,\prod_v L(s,\pi_v\times \tau_v)$$
and
$$\varepsilon(s,\pi\times\tau)=\,\prod_v \varepsilon(s,\pi_v\times\tau_v,\psi_v).$$
The $L$-function $L(s,\pi\times\tau)$ is said to be {\bf nice} if it satisfies the following properties (N1)-(N-3):

(N1) $L(s,\pi\times\tau)$ and $L(s,\pi^\vee\times\tau^\vee)$ has A.C. to the whole complex plane;
\\
(N2) $L(s,\pi\times\tau)$ and $L(s,\pi^\vee\times\tau^\vee)$ are bounded in vertical strips of finite width;
\\
(N3) These entire functions satisfy the standard functional equation
\begin{center}
\framebox[15cm][c]{$L(s,\pi\times\tau)=\,\varepsilon(s,\pi\times\tau)\,L(1-s,\pi^\vee\times\tau^\vee)$}
\end{center}

\end{slide}



\newpage
\begin{slide}

{\large\bf Converse Theorem (Cogdell and\\ Piatetski-Shapiro [1994,\,1999]).}
\par
Let $\pi$ be an irred admissible repn of $GL(n)$ whose central character is trivial on $F^*$ and whose $L$-function $L(s,\pi)$ converges absolutely in some half plane. Let
$S$ be a finite set of finite places. Assume that $L(s,\pi\times \tau)$ is {\bf nice} for every cuspidal repn $\tau$ of $GL(m)$ for $1\leq m\leq n-2$, which is unramified at the places in $S$. Then there is an auto repn $\pi'$ of $GL(n)$ such that $\pi_v\cong \pi_v'$ for all $v\notin S.$

\vskip 1cm
\noindent {\large\bf Example 3.} (Tensor Product)
\par
(1) $GL(2)\times GL(2)\lrt GL(4)$\\ \indent \ \ \ \ \ (D. Ramakrishnan [2000])\\
(2) $GL(2)\times GL(3)\lrt GL(6)$ \\ \indent \ \ \ \ \ \ \ (Kim-Shahidi [2002]).
\end{slide}


\begin{slide}
\newpage

\smallskip
\noindent {\large\bf Example 4.} (Asgari-Shahidi [2005])

$H=GSpin(m),\ G=GL(N)$ with $N=m$ or $2[ {\frac m2}]$ and $\sigma$ is an embedding. We note that
$$\LH^0=
\begin{cases} GSO(m) & \quad {\rm if}\ m\ {\rm is\ even}\\
GSp\big(2[ {\frac m2}]\big)& \quad {\rm if}\ m \ {\rm is\ odd.}
\end{cases}$$
${\bf LF}$ holds for an irreducible {\it generic} cuspidal repn of $H(\A)$.

\vskip 1cm
\noindent {\large\bf Example 5.} (Jiang-Soudry [2003])

\par
If $H=SO(2n+1)$ and $G=GL(2n)$,\\
the functorial lift from $H$ to $G$ is {\bf injective} up to isomorphism.
\end{slide}


\begin{slide}
\newpage

\noindent {\large\bf Example 6.} (Symmetric Power Product for $GL(2)$).\\
Let $H=GL(2)$ and $G=GL(m+1).$ Then we have an $L$-homomorphism
$${\rm Sym}^m: GL(2,\BC)\lrt GL(m+1,\BC).$$
(1) $m=2$ (Gelbart-Jacquet [1978]).\\
(2) $m=3$ (Kim-Shahidi [2002])\\
(3) $m=4$ (Kim [2003]).\par\noindent
It was observed by Langlands that the functoriality for ${\rm Sym}^m$ for all $m\geq 1$ implies the {\bf Ramanjujan Conjecture} for Maass forms (also the {\bf Selberg Conjecture}) and the {\bf Sato-Tate Conjecture}.
\end{slide}

\newpage
\begin{slide}
\noindent {\large\bf Example 7.} (Exterior Square for $GL(4)$)
$$GL(4)\longrightarrow GL(6)$$
\ \ \ \ \ \ \ \ \ \ \ \ \ \ \ \ \ \ \ (Kim [2003])
\par\noindent {\large\bf Example 8.} (P.-S. Chan and Y. Flicker [2007])\par\noindent
$\bullet \ E/F$ : a quadratic ext of a number field $F$\\
$\bullet \ F'/F$ : a cyclic ext of an odd degree\\
Then there is a {\bf BC\ functorial\ lifting} from $U(3,E/F)$ to $U(3,F'E/F').$\par\noindent $\large\verb"Tools of Proof:"$\\
$\circledcirc$ Trace Formula\\
$\circledcirc$ Known results on BC lifiting from $U(3,E/F)$ to $GL(3,E)$ and on BC lifting for $GL(n)$.
\end{slide}


\newpage
\begin{slide}
\begin{center}
{\large \bf 3. Applications of Functoriality}
\end{center}

\noindent {\bf [1] The Ramanujan Conjecture}\par\noindent
(1-1) {\bf R.C for holomorphic cusp forms}\\
Let
$$f(\tau)=\sum_{n\geq 1} a_n\,e^{2\pi in\tau}\ (a_1=1)$$ be a holomorphic cusp new eigenform for $\Gamma_0(N)$. Then
by P. Deligne [1974], R.C. is valid for $f$, in other words, for $p\nmid N$,
$$|a_p|\leq 2\,p^{{k-1}\over 2}.$$

(1-2) {\bf R.C. for Maass cusp forms}\\
Let
$$f(x+iy)=\sum_{n\neq 0} a_n\,K_{s/2}(2\pi |n|y)(|n|y)^{1/2}\,e^{2\pi inx}$$
be the Fourier expansion of a Maass cusp form  $f$ for $\Gamma_0(N)$ such that
$$\Delta f ={\frac 14} \big(1-s^2\big)f, \quad s\in\BC,$$
\end{slide}

\begin{slide}
where $$\Delta=-y^2 \big(\partial_x^2+\partial_y^2\big)\ \ {\rm Laplacian}.$$
For $p\nmid N$, we write
$$a_p=p^{-1/2}\big( \alpha_p+\alpha_p^{-1}\big),\quad \alpha_p\in\BC.$$\par
\noindent
\begin{center}
\framebox[16cm][c]{{\bf Ramanujan Conj\,:}\ \ $|\alpha_p|=1$ for all $p\nmid N.$}
\end{center}
\par\noindent
{\bf Theorem (Kim-Sarnak [2003]).} Let $N=1.$ For any prime $p$,
\begin{equation} p^{-{7\over {64}}}\leq |\alpha_p| \leq p^{7\over {64}}.\end{equation}
It is the best known estimate so far. The above theorem is a consequence of some results in Langlands functoriality.
By the estimate (2), the eigenvalue $\lambda$ of a Maass form $f$ is given by
\begin{equation}
\lambda\geq {\frac 14}- \left( {\frac 7{64}}\right)^2\approx 0.238037\cdots.
\end{equation}
\end{slide}

\newpage
\begin{slide}
\noindent {\bf [2] The Sato-Tate Conjecture}\\

Suppose $f$ is a a holomorphic cusp new eigenform for $\Gamma_0(N)$ as above \\
or\\
Suppose $f$ is a Maass cusp form with the assumption of the validity of the Ramanujan Conjecture.\\
\ \ We set
$$ \alpha_p = e^{i \theta_p},\quad \theta_p\in [0,\pi].$$

\noindent {\bf Sato-Tate Conjecture.} The sequence $\big\{ \theta_p \big\}$ is equi-distributed (or uniformly distributed) with respect to the measure


\begin{center}
\framebox[5cm][c]{${2\over {\pi}}\,\sin^2 t \,dt$}
\end{center}

on $[0,\pi].$ In other words, for any real numbers $\alpha$ and $\beta$ with $0\leq \alpha < \beta\leq \pi,$
$$\int_\alpha^{\beta} {2\over {\pi}}\,\sin^2 t \,dt=\lim_{X\longrightarrow \infty}
{ {|\{ p\leq X\,|\ p\nmid N,\ \alpha < \theta_p < \beta \}|} \over {|\{ p\,|\ p\leq X,\ p\nmid N\,\}|} }.$$
\end{slide}

\begin{slide}
Recently L. Clozel, M. Harris and R. Taylor [2006] proved the Sato-Tate conjecture for elliptic curves with a certain condition.
\par\noindent
{\bf Theorem [Clozel-Harris-Taylor].} Let $E/\BQ$ be an elliptic curve of non-CM type. Assume $j(E)\notin \BZ.$ Then the Sato-Tate conjecture is valid for $E$.\par\noindent
{\bf Ideas of Proof.} By the modularity of $E$, the $L$-function $L(s,E)$ is entire. For a prime $p$, we put
$$a_p=p+1- | E({\mathbb F}_p)|.$$
Then
$$ L(s,E)=\prod_p L_p(s,E),$$
where
$$L_p(s,E)=\begin{cases} \big( 1-a_p p^{-s}\big)^{-1} & {\rm if}\ p\ {\rm is\ bad,}\\
                          \big( 1-a_p p^{-s}+p^{1-2s}\big)^{-1} & {\rm if}\ p\ {\rm is\ good.}
                          \end{cases}$$
\end{slide}

\begin{slide}
For a good prime $p$, Hasse's inequality gives
$$|a_p|\leq 2\, p^{1/2}.$$
We set
\begin{center}
\framebox[8cm][c]{$a_p^*= {{a_p}\over {2\,p^{1/2}} }\in [-1,1]$}
\end{center}
We can also write
$$a_p^*=\cos\theta_p\quad {\rm for\ a\ unique}\ \theta_p\in [0,\pi].$$
In this special case (the holomorphic cusp eigenform of weight 2 corresponding to $E$), the Sato-Tate conjecture can be rephrased as follows:
\par\noindent
{\bf Sato-Tate Conjecture.} Suppose $E/\BQ$ has no complex multiplication. Then the sequence $\{ a_p^*\}$ (resp. $\{\theta_p\}$)
is equidistributed in [-1,1] (resp. $[0,\pi]$) with respect to the probability measure
\begin{center}
\framebox[12cm][c]{${2\over {\pi}}\,\sqrt{1-t^2}\,dt\quad \left({\rm resp.}\ {2\over {\pi}}\,\sin^2\theta\,d\theta\,\right)$}
\end{center}

\end{slide}

\begin{slide}
Let $\alpha_p$ and $\beta_p$ be complex numbers such that
$$\alpha_p+\beta_p=a_p\quad {\rm and}\quad    \alpha_p\,\beta_p=p.$$
For a good prime $p$, we see that
$$L_p(s,E)=\Big[ \big( 1-\alpha_p\, p^{-s}\big) \big( 1-\beta_p\, p^{-s}\big) \big]^{-1}.$$
For $n\in \BZ^+$ and a good prime $p$, we put
$$L_p(s,E,{\rm Sym}^n)=\left[ \prod_{i=0}^n \big( 1-\alpha_p^i\,\beta_p^{n-i}\,p^{-s}\big)\right]^{-1}.$$
There is also a definition of $L_p(s,E,{\rm Sym}^n)$ for a bad prime $p$.
\par\noindent
Define
$$ L(s,E,{\rm Sym}^n) =\prod_p L_p(s,E,{\rm Sym}^n).$$
If $n=1,\ L(s,E,{\rm Sym})$ is entire and nonvanishing on ${\rm Re}\,( s)=1$ (R. Rankin).\par\noindent
If $n=2,\ L\left(s,E,{\rm Sym}^2\right)$ is entire (G. Shimura [$F=\BQ$ : 1975], Gelbart-Jacquet [$F$=any number field: 1978]).
\end{slide}

\begin{slide}
\par\noindent
If $n=3,\ L\left(s,E,{\rm Sym}^3\right)$ is entire (Kim-Shahidi [2002]).
\par\noindent
If $n=4,\ L\left(s,E,{\rm Sym}^4\right)$ is entire (Kim [2003]).
\par\noindent
If $n=5,6,7,8,\ L\left(s,E,{\rm Sym}^n\right)$ extends to a meromorphic function to $\BC$ which is holomorphic in
${\rm Re}\,( s)\geq 1$ (Kim-Shahidi [2002]).
\par\noindent
If $n=9,\ L\left(s,E,{\rm Sym}^9\right)$ extends to a meromorphic function to $\BC$ which is holomorphic in
${\rm Re}\,( s)> 1$ and may have a pole at $s=1$ (Kim-Shahidi [2002]).
\par\noindent
{\bf Theorem\,(J.-P. Serre).} Let $E/\BQ$ be an elliptic curve over $\BQ$. Assume for any positive integer $n\in\BZ^+,\
L\left(s,E,{\rm Sym}^n\right)$ extends to a meromorphic function on $\BC$ which is holomorphic and nonvanishing in
${\rm Re}\,( s)\geq 1+ {\frac n2}.$ Then the Sato-Tate conjecture is valid for $E.$
\end{slide}

\newpage
\begin{slide}
When $j(E)\notin \BZ$, Clozel, Harris and Taylor proved that for any positive integer $n\in\BZ^+,\
L\left(s,E,{\rm Sym}^n\right)$ extends to a meromorphic function on $\BC$ which is holomorphic and non-vanishing in
${\rm Re}\,( s)\geq 1+ {\frac n2}.$ By Serre's Theorem, the Sato-Tate conjecture is true for $E$.


In order to prove this fact, when $n$ is {\bf odd}, they first proved that $L\left(s,E,{\rm Sym}^n\right)$ is {\bf potentially automorphic}, in other words,
$L\left(s,E,{\rm Sym}^n\right)$ is associated to a cuspidal rep of $GL(n+1)$ over some totally real Galois ext of $\BQ$.
This argument involves two parts.
\par\noindent
(1) An extension of Wiles' technique used to prove the modularity of an elliptic curve based on Galois cohomology and analysis of automorphic representations on different groups, in particular unitary groups.
\end{slide}

\newpage
\begin{slide}
(2) An extension of Taylor's idea for proving meromorphic continuation of $L$-functions attached to 2-dimensional Galois representations.
\par
For any {\bf even} $n$, Harris, Shepherd and Taylor found a twisted form of the moduli space of certain Calabi-Yau manifolds that could be used to study $n$-dimensional Galois representations.
\hfill $\square$
\par\noindent
{\bf Theorem [Taylor,\ 2008].} Let K be a totally real field and let $E/K$ an elliptic curve with multiplicative reduction at some prime. Then the Sato-Tate conjecture is {\bf true} for $E$, i.e., the numbers
$${ {1+p-|E({\mathbb F}_p)|}\over {2p^{1/2}} }$$
are equidistributed in $[-1,1]$ w.r.t the probability measure
$$ {2\over {\pi}}\,\sqrt{1-t^2}\,dt.$$
\end{slide}

\newpage
\begin{slide}
\begin{center}
{\large \bf References }
\end{center}

\vskip 0.51cm

[1] $ \textbf{L. Clozel, M. Harris and R. Taylor}$, {\em  Automorphy for some $\ell$-adic lifts of
automorphic mod $\ell$ Galois representations}, Publ. Math. IHES {\bf 108} (2008), 1-181.

[2] $\textbf{R. Taylor}$, {\em  Automorphy for some $\ell$-adic lifts of
automorphic mod $\ell$ Galois representations. II},
Publ. Math. IHES {\bf 108} (2008), 183-239.

[3] $ \textbf{M. Harris, N. Shepherd-Barron and}$\\ $\textbf{R. Taylor}$,
{\em A family of Calabi-Yau varieties and potential automorphy}, to appear in
Annals of Math.\,(2009).

[4] $ \textbf{M. Harris, N. Shepherd-Barron and}$\\ $\textbf{R. Taylor}$,
{\em Ihara's lemma and potential automorphy},
preprint (2006).

\end{slide}

\newpage
\begin{slide}

[5] $ \textbf{M. Harris and R. Taylor}$,
{\em The geometry and cohomology of some simple Shimura varieties},
Annals of Math. Studies {\bf 151}, Princeton Univ. Press (2001).

[6] $\textbf{M. Kisin}$, {\em  Moduli of finite flat group schemes, and modularity}, to appear
in Annals of Math.

[7] $ \textbf{C. Skinner and A. Wiles}$,
{\em Base Change and a problem of Serre},
Duke Math. Journal {\bf 107} (2001), 15-25.

[8] $ \textbf{A. Wiles}$,
{\em Modular elliptic curves and Fermat's last theorem},
Annals of Math. {\bf 141} (1995), 443-551.

\end{slide}

\newpage
\begin{slide}
\noindent {\large\bf [3] The Inverse Galois Problem}\\

\ \ Using
the functorial lifting from $SO(2n+1)$ to $GL(2n)$, C.Khare, M. Larsen
and G. Savin [2008] proved that for any prime $\ell$ and any
even positive integer $n$, there are infinitely many exponents $k$
for which the finite simple group $PSp_n({\mathbb F}_{\ell^k})$
appears as a Galois group over $\BQ$. \par\noindent
$\bullet$ Functoriality and the inverse Galois problem, Compositio Math. {\bf 144} (2008), 541-564.
\\

Furthermore, in their recent
paper [2008] they extended their earlier work to prove that
for a positive integer $t$, assuming that $t$ is even if $\ell=3$
in the first case (1) below, the following statements (1)-(3)
hold: \vskip 0.15cm (1) Let $\ell$ be a prime. Then there exists
an integer $k$ divisible by $t$ such that the simple group
$G_2({\mathbb F}_{\ell^k})$ appears as a Galois group over
$\BQ$.
\end{slide}

\newpage
\begin{slide}
(2) Let $\ell$ be an odd prime. Then there exists an
integer $k$ divisible by $t$ such that the simple finite group
$SO_{2n+1}({\mathbb F}_{\ell^k})^{\rm der}$  or the finite
classical group $SO_{2n+1}({\mathbb F}_{\ell^k})$ appears as a
Galois group over
$\BQ$.\par
\indent (3) If $\ell\equiv 3,5\, ({\rm mod}\,8)$ and $\ell$ is a
prime, then there exists an integer $k$ divisible by $t$ such that
the simple finite group $SO_{2n+1}({\mathbb F}_{\ell^k})^{\rm
der}$ appears as a Galois group over $\BQ$.

\vskip 0.2cm The construction of Galois groups in (1)-(3) is based
on the functorial lift from $Sp(2n)$ to $GL(2n+1)$, and the
backward lift from $GL(2n+1)$ to $Sp(2n)$ plus the theta lift from
$G_2$ to $Sp(6)$.
\par\noindent
$\bullet$ Functoriality and the inverse Galois problem {\rm II}:
Groups of type $B_n$ and $G_2$, arXiv:0807.0861v1 (2008).
\end{slide}

\newpage
\begin{slide}
\noindent {\large\bf [4] Cuspidality and Irreducibility}
\par\noindent
{\bf Question:} When is a functorial lift $\Pi$ of $\pi$ {\it cuspidal} ? Give a criterion that $\Pi$ is {\it cuspidal}.
\par We have partial results for the special cases.

\par\noindent
{\large\bf Example 1.} Let $\pi$ be a cuspidal rep of $GL(2).$ Then the following properties are known:\\
(1) $\pi$ is dihedral iff ${\rm Sym}^2(\pi)$ is not cuspidal. (G-J [1978]).\\
(2) Assume that $\pi$ is not dihedral and not monomial. Then $\pi$ is tetrahedral if and only if
${\rm Sym}^3(\pi)$ is not cuspidal. (K-S [2002]).\\
(3) Assume that $\pi$ is not dihedral and not monomial. Then $\pi$ is not octahedral if and only if
${\rm Sym}^4(\pi)$ is cuspidal. (Kim [2003]).\par
\indent By (2) and (3), $\pi$ is octahedral if and only if ${\rm Sym}^3(\pi)$ is cuspidal but ${\rm Sym}^4(\pi)$ is not cuspidal.
\end{slide}

\newpage
\begin{slide}
{\large {\bf Definitions:}} (1) $\pi$ is called {\bf monomial} if $\pi\simeq \pi \otimes \eta$
for some idele class character $\eta$ of $F$.
\par
(2) $\pi$ is said to be {\bf solvable\ polyhedral} if it is dihedral, tetrahedral or octahedral.
\par
(3) $\pi$ is said to be {\bf self-dual} if $\pi\simeq \pi^{\vee}.$
\par
(4) $\pi$ is said to be {\bf essentially\ self-dual} if $\pi^{\vee}\simeq \pi\otimes \chi$
for some idele class character $\chi$.
\par
(5) $\pi$ is said to be {\bf almost\ self-dual} if $\pi^{\vee}\simeq \pi\otimes \xi\,|\,\cdot\,|^t$, where $\xi$ is a finite order character and $t\in\BC.$
\par
(6) According to Langlands, the archimedean component $\pi_\infty$ is associated to an $n$-dimensional representation $\sigma_\infty(\pi)$
of the real Weil group $W_\BR$. We have a canonical exact sequence
\begin{equation*}
1\lrt \BC^* \lrt W_\BR \lrt {\textrm{Gal}}(\BC/\BR)\lrt 1
\end{equation*}

\end{slide}


\newpage
\begin{slide}

which represents the unique nontrivial extension of ${\textrm{Gal}}(\BC/\BR)$ by $\BC^*.$ One has a decomposition
\begin{equation}
\sigma_\infty (\pi)|_{\BC^*}\simeq \,\oplus_{j=1}^n\chi_j,
\end{equation}
where each $\chi_j$ is a quasi-character of $\BC^*.$ $\pi$ is said to be {\bf algebraic} if each $\chi_j$ is
algebraic, that is, there are integers $p_j,q_j$ such that
$$\chi_j(z)=\,z^{p_j}\,{\overline z}^{q_j},\quad z\in\BC^*.$$
We say that $\pi$ is {\bf regular} if the multiplicity of each $\chi_j$ is one, and that $\pi$ is
{\bf semi-regular} if the multiplicity of each $\chi_j$ is at most two.
\par
(7) An isobaric automorphic representation $\pi$ of $GL(n,{\mathbb A})$ is called {\bf odd} if for every one-dimen'l representation $\xi$
of $W_\BR$ occurring in $\sigma_\infty (\pi)$, the following condition
$$|m_+(\pi,\xi)-m_-(\pi,\xi)|\leq 1,$$
where $m_+(\pi,\xi)$ (resp. $m_-(\pi,\xi)$) denotes the multiplicity of the eigenvalue $+1$ (resp. $-1$) under the action
$\BR^*/\BR^*_+$ on the $\xi$-isotypic component of $\sigma_\infty (\pi)$.

\end{slide}

\newpage
\begin{slide}
\noindent
{\large\bf Example 2} (D. Ramakrishnan and S. Wang).

 Let $\pi$ and $\theta$ be cuspidal reps of $GL(2)$ and $GL(3)$ respectively over a number field $F$. Then the functorial lift $\pi\boxtimes \theta$ on $GL(6)/F$ is {\bf cuspidal} unless one of the following satisfies
\par\noindent (a) $\pi$ is dihedral, and $\theta$ is a twist of $Ad(\pi);$
\par\noindent (b) $\pi$ is dihedral, $L(s,\pi)=L(s,\chi)$ for an idele class character $\chi$ of a cubic, non-normal ext $K$ of $F$, and the BC $\pi_K$ is Eisensteinian.
\par\noindent
{\bf Remark.} The $L$-function $L(s,\pi\boxtimes \theta)$ is equal to the Rankin-Selberg $L$-function $L(s,\pi\times\theta).$ By Example 1, $\pi\boxtimes\theta$ is {\it cuspidal} if $\pi$ is not dihedral and $\theta$ is not a twist of $Ad(\pi).$
\end{slide}

\newpage
\begin{slide}
\par\noindent
{\large\bf Example 3}\,(D. Ramakrishnan and S. Wang).

Let $\pi$ and $\theta$ be cuspidal reps of $GL(2)$ and $GL(3)$ respectively over a number field $F$. Assume that $\pi$ is not of solvable polyhedral type, and $\theta$ not essentially self-dual. Then we have the following

(a) If $\theta$ does not admit any self twist, the functorial lift $\pi\boxtimes\theta$ is {\it cuspidal} without any self twist.

(b) If $\theta$ is not of solvable type, $\pi\boxtimes\theta$ is {\it cuspidal} and not of solvable type.
\end{slide}

\newpage
\begin{slide}
\noindent
{\large\bf Example 4}\,(D. Ramakrishnan [2007]).

Let $\pi$ be a cuspidal rep of $GL(2)/F$ with central character $\omega$. Assume that $\pi$ is not solvable polyhedral and also that $\textrm{Sym}^m(\pi)$ is automorphic for all $m\geq 1.$ Then

(a) $\textrm{Sym}^5(\pi)$ is {\it cuspidal};

(b) $\textrm{Sym}^6(\pi)$ is {\it not\ cuspidal} if and only if we have
$$ \textrm{Sym}^5(\pi) \cong Ad(\tau)\boxtimes \pi \otimes \omega^{-4}$$
for a cuspidal rep $\tau$ of $GL(2)/F$;

(c) If $\textrm{Sym}^6(\pi)$ is cuspidal, then so is $\textrm{Sym}^m(\pi)$ for all $m\geq 1$;

(d) If $F=\BQ$ and $\pi$ is a holomorphic, non-CM new eigenform of weight $k\geq 2$, then
$\textrm{Sym}^m(\pi)$ is cuspidal for all $m\geq 1$.

\end{slide}

\newpage
\begin{slide}
\noindent
{\large\bf Example 5}\,(D. Ramakrishnan [2008]).

Let $n\leq 5$ and $\ell$ be a prime. Suppose a continuous $\ell$-adic representation $\rho_{\ell}$ of the absolute Galois group $G_\BQ$ is
associated to an isobaric automorphic representation $\pi$ of $GL(n,{\mathbb A})$.  Assume the following conditions (a),(b) and (c):
\par
\indent (a) $\rho_{\ell}$ is irreducible.
\par
\indent (b) $\pi$ is odd if $n\geq 3$.
\par
\indent (c) $\pi$ is semi-regular if $n=4$, and regular if $n=5.$
\par
Then $\pi$ is {\bf cuspidal}.

\end{slide}

\newpage
\begin{slide}
\begin{center}
{\large \bf 4. Methods to Tackle LFC}
\end{center}

\par\noindent
{\bf $\circledcirc$ Trace Formula}

\indent\ \ \ [J. Arthur, R. Langlands,$\cdots$]

\par\noindent
{\bf $\circledcirc$ Converse Theorem}

\indent\ \ \ [Cogdell, Piateski-Shapiro\,(1929-2009),$\cdots$]

\par\noindent
{\bf $\circledcirc$ Endoscopy}

\indent\ \ \ [R. Langlands, D. Shelstad, R. Kottwitz,$\cdots$]

\par\noindent
{\bf $\circledcirc$ Base Change}

\indent\ \ [R. Langlands, J. Arthur, L. Clozel, Y. Flicker,$\cdots$]

\par\noindent
{\bf $\circledcirc$ Theta Correspondence}

\indent\ \ \ [R. Howe, J.-S. Li,$\cdots$]

\par\noindent
{\bf $\circledcirc$  Motif and Shimura Varieties}

\end{slide}


\newpage
\begin{slide}

\par\noindent
$\circledcirc$ {\Large\bf New Method (Idea)\ ???}

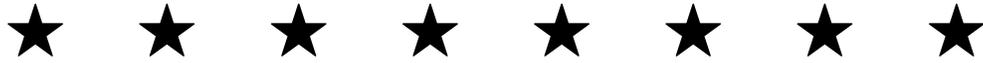
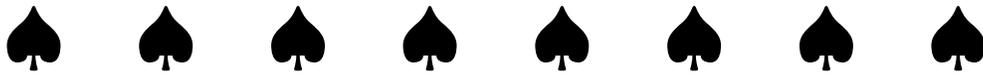
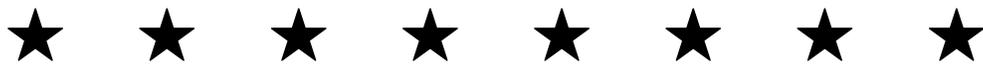
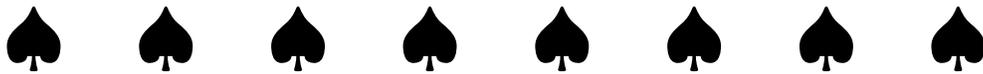
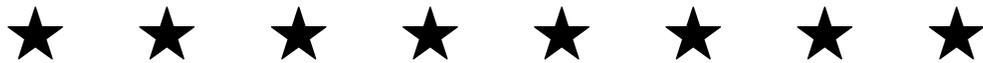
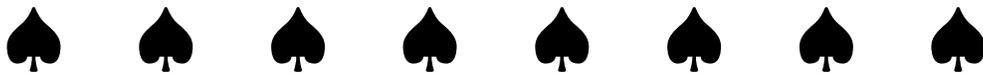
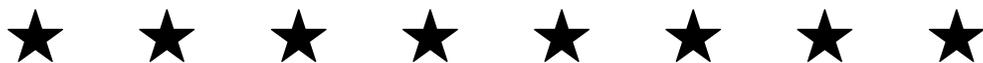
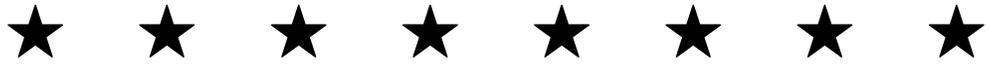
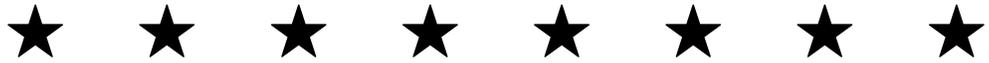
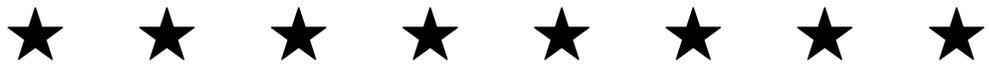
\begin{center}
{\Large $ \textbf{$\bigstar\ \ \bigstar\ \ \bigstar\ \ \bigstar\ \ \bigstar
\ \ \bigstar\ \ \bigstar\ \ \bigstar$}$}
\end{center}

\begin{center}
{\Large $ \textbf{$\spadesuit\ \ \spadesuit\ \ \spadesuit\ \ \spadesuit\ \ \spadesuit
\ \ \spadesuit\ \ \spadesuit\ \ \spadesuit$}$}
\end{center}

\begin{center}
{\Large $ \textbf{$\bigstar\ \ \bigstar\ \ \bigstar\ \ \bigstar\ \ \bigstar
\ \ \bigstar\ \ \bigstar\ \ \bigstar$}$}
\end{center}

\begin{center}
{\Large $ \textbf{$\spadesuit\ \ \spadesuit\ \ \spadesuit\ \ \spadesuit\ \ \spadesuit
\ \ \spadesuit\ \ \spadesuit\ \ \spadesuit$}$}
\end{center}

\begin{center}
{\Large $ \textbf{$\bigstar\ \ \bigstar\ \ \bigstar\ \ \bigstar\ \ \bigstar
\ \ \bigstar\ \ \bigstar\ \ \bigstar$}$}
\end{center}

\begin{center}
{\Large $ \textbf{$\spadesuit\ \ \spadesuit\ \ \spadesuit\ \ \spadesuit\ \ \spadesuit
\ \ \spadesuit\ \ \spadesuit\ \ \spadesuit$}$}
\end{center}

\begin{center}
{\Large $ \textbf{$\bigstar\ \ \bigstar\ \ \bigstar\ \ \bigstar\ \ \bigstar
\ \ \bigstar\ \ \bigstar\ \ \bigstar$}$}
\end{center}


\end{slide}


\newpage

\begin{slide}

\begin{center}
{\Large $ \textbf{$\bigstar\ \ \bigstar\ \ \bigstar\ \ \bigstar\ \ \bigstar
\ \ \bigstar\ \ \bigstar\ \ \bigstar$}$}
\end{center}

\begin{center}
{\Large $ \textbf{$\bigstar\ \ \bigstar\ \ \bigstar\ \ \bigstar\ \ \bigstar
\ \ \bigstar\ \ \bigstar\ \ \bigstar$}$}
\end{center}

\begin{center}
{\Large $ \textbf{$\bigstar\ \ \bigstar\ \ \bigstar\ \ \bigstar\ \ \bigstar
\ \ \bigstar\ \ \bigstar\ \ \bigstar$}$}
\end{center}

\vskip 1cm

\begin{center}  \textbf{\Large Thank You Very Much !!!}\end{center}

\vskip 1cm

\begin{center}
{\Large $ \textbf{$\bigstar\ \ \bigstar\ \ \bigstar\ \ \bigstar\ \ \bigstar
\ \ \bigstar\ \ \bigstar\ \ \bigstar$}$}
\end{center}

\begin{center}
{\Large $ \textbf{$\bigstar\ \ \bigstar\ \ \bigstar\ \ \bigstar\ \ \bigstar
\ \ \bigstar\ \ \bigstar\ \ \bigstar$}$}
\end{center}

\begin{center}
{\Large $ \textbf{$\bigstar\ \ \bigstar\ \ \bigstar\ \ \bigstar\ \ \bigstar
\ \ \bigstar\ \ \bigstar\ \ \bigstar$}$}
\end{center}

\begin{center}
{\Large $ \textbf{$\bigstar\ \ \bigstar\ \ \bigstar\ \ \bigstar\ \ \bigstar
\ \ \bigstar\ \ \bigstar\ \ \bigstar$}$}
\end{center}

\end{slide}

\end{document}